\documentclass[pdflatex,sn-mathphys-num]{sn-jnl}

\usepackage{graphicx}
\usepackage{amsmath,amssymb,amsfonts,amsthm}
\usepackage{algorithm,algorithmic}
\usepackage{multirow,booktabs,array,longtable,enumitem}
\usepackage{xcolor,microtype,placeins,needspace}
\usepackage{aliascnt}

\theoremstyle{thmstyleone}
\newtheorem{theorem}{Theorem}[section]
\newtheorem{lemma}{Lemma}[section]
\newtheorem{proposition}{Proposition}[section]

\newaliascnt{assumption}{theorem}
\newtheorem{assumption}[assumption]{Assumption}
\aliascntresetthe{assumption}
\theoremstyle{thmstyletwo}
\newtheorem{remark}{Remark}[section]
\theoremstyle{thmstylethree}
\newtheorem{definition}{Definition}[section]

\usepackage[capitalise,noabbrev]{cleveref}
\crefname{assumption}{Assumption}{Assumptions}
\Crefname{assumption}{Assumption}{Assumptions}
\crefname{theorem}{Theorem}{Theorems}
\Crefname{theorem}{Theorem}{Theorems}
\crefname{lemma}{Lemma}{Lemmas}
\Crefname{lemma}{Lemma}{Lemmas}
\crefname{proposition}{Proposition}{Propositions}
\Crefname{proposition}{Proposition}{Propositions}
\crefname{corollary}{Corollary}{Corollaries}
\Crefname{corollary}{Corollary}{Corollaries}
\crefname{definition}{Definition}{Definitions}
\Crefname{definition}{Definition}{Definitions}
\crefname{remark}{Remark}{Remarks}
\Crefname{remark}{Remark}{Remarks}
\crefname{example}{Example}{Examples}
\Crefname{example}{Example}{Examples}

\hypersetup{hypertexnames=false,colorlinks=true,linkcolor=blue,citecolor=blue,urlcolor=blue,
  bookmarksopen=true,bookmarksnumbered=true}

\newcommand{\R}{\mathbb{R}}

\newcommand{\Xstar}{X^{*}}
\newcommand{\dom}{\operatorname{dom}}
\newcommand{\ri}{\operatorname{ri}}
\newcommand{\argmin}{\operatorname*{arg\,min}}
\newcommand{\ip}[2]{\langle #1,#2\rangle}
\newcommand{\norm}[1]{\lVert #1\rVert}

\newcommand{\Acrit}{\mathcal{A}}
\newcommand{\Ncorr}{\mathcal{N}}
\newcommand{\EE}{\mathcal{E}}
\newcommand{\HH}{\mathcal{H}}
\newcommand{\hh}{\mathfrak{h}}
\newcommand{\bb}{\mathfrak{b}}

\newcolumntype{P}[1]{>{\raggedright\arraybackslash}p{#1}}
\begin{document}
\title[General Inertia and Critical Rates]{The Convergence Rate of Accelerated Proximal Gradient Methods with General Inertial Term for Convex Composite Multiobjective Optimization}

\author[1]{\fnm{Chengzhi} \sur{Huang}}\email{huangczmaths@163.com}

\affil[1]{\orgname{Sichuan University}, \orgaddress{\city{Chengdu}, \postcode{610065}, \state{Sichuan}, \country{China}}}

\abstract{
We study accelerated proximal gradient methods with general stepsizes for convex composite multiobjective optimization, where each objective is a smooth convex term plus a proper closed convex term and the proximal subproblem keeps a separate offset for each objective. For a two-parameter root momentum rule with trajectory-dependent coefficients, we construct a new energy function and establish an energy identity. This yields sharper convergence guarantees: boundedness, weighted square summability of the displacements, a merit-function rate of $o(k^{-2})$, and displacement decay of $o(k^{-1})$ whenever the first parameter converges to a value in $(0,1)$; the second parameter need not converge or have bounded total variation. At the critical value $\alpha=3$, we obtain
$u_0(x^N)\le \frac{2\ell R}{N(N+3)},$
together with $O(1/N)$ displacement decay, square summability, and weak Pareto convergence. The remaining little-$o$ question at $\alpha=3$ is reduced to the vanishing of a nonnegative critical energy. For $1<\alpha<3$, we prove the rate $O(N^{-(\alpha-1)})$, and for every $\alpha>0$ we obtain $u_0(x^k)\to0$, $F_i(x^k)\to c_i$, and $\sum_k D_k/k<\infty$. A JOS1 experiment numerically validates the predicted little-$o$ separation and the rate behavior.
}

\keywords{composite multiobjective optimization, accelerated proximal gradient method, momentum factor, general stepsizes, critical parameter, subcritical rate}

\pacs[MSC Classification]{49J52, 65K05, 90C25, 90C29}

\maketitle

\section{Introduction}\label{sec:intro}

We consider a composite multiobjective optimization problem of the form
\begin{equation}\label{eq:intro-problem}
 \min_{x\in\R^n}F(x)=\bigl(F_1(x),\ldots,F_m(x)\bigr),
 \qquad F_i=f_i+g_i,\qquad i=1,\ldots,m.
\end{equation}
Each $f_i:\R^n\to\R$ is convex and continuously differentiable with a globally Lipschitz gradient, while each $g_i:\R^n\to\R\cup\{+\infty\}$ is proper, closed, and convex. A natural way to solve this problem is via a proximal gradient method: in each multiobjective proximal gradient step, every $f_i$ is linearized at a common point, each $g_i$ is kept inside the subproblem, and the resulting componentwise models are combined by taking their maximum \cite{tanabe2019proximal}. This construction does not prescribe scalarization weights in advance and retains the computational structure of the $g_i$. To achieve faster solving, acceleration techniques are introduced. Although acceleration often improves practical performance, the convergence rates established for existing multiobjective accelerated algorithms do not always match their observed numerical performance. This discrepancy motivates the vector-specific design question addressed next.

A central design question is therefore: how much freedom may the momentum coefficient have while the accelerated rate for the original merit function is preserved? The scalar case offers a complete answer in terms of a single parameter $\alpha$. The inertial system with asymptotic vanishing damping
\begin{equation}\label{eq:intro-avd}
 \ddot x(t)+\frac{\alpha}{t}\dot x(t)+\nabla\Phi(x(t))=0
\end{equation}
has the rate $O(t^{-2\alpha/3})$ for $\alpha\le3$, the value $O(t^{-2})$ at $\alpha=3$, and $o(t^{-2})$ for $\alpha>3$; the critical value $\alpha=3$ therefore separates a plateau of fast rates from a linearly increasing subcritical regime \cite{attouch2017subcritical}. For composite multiobjective problems the corresponding picture is far from complete, and it is the purpose of this paper to supply it for exact proximal steps and general stepsizes.


For a single objective, the proximal gradient method combines a gradient step on the smooth part with a proximal step on the nonsmooth part, and FISTA adds the Nesterov momentum factor $t_{k+1}=(1+\sqrt{1+4t_k^2})/2$ to obtain the optimal worst-case rate $O(1/k^2)$ \cite{beck2009fista}. The multiobjective extension replaces the scalar objective value by a componentwise model and solves one offset subproblem per iteration; its basic version converges with rate $O(1/k)$ \cite{tanabe2019proximal}, and its accelerated version recovers $O(1/k^2)$ \cite{tanabe2023accelerated,tanabe2023convergence}.

Two later developments motivate the present analysis. First, the momentum factor itself can be generalized: for constants $a\in[0,1)$ and $b\in[a^2/4,1/4]$, the root recurrence $t_{k+1}=\sqrt{t_k^2-at_k+b}+1/2$ covers the classical rule and preserves the $O(1/k^2)$ rate \cite{tanabe2022root}. Second, the rational coefficient $\beta_k=(k-1)/(k+\alpha-1)$ with $\alpha>3$ leads to an $o(k^{-2})$ rate for multiobjective proximal gradient methods \cite{zhang2023}. In the inertial-gradient literature the same coefficient is the discrete counterpart of the vanishing damping in \eqref{eq:intro-avd}, and the multiobjective algorithms obtained from inertial gradient-like systems are analysed through Lyapunov functions \cite{sonntag2024fast}.

These results concern parameters that are fixed before the computation, or that follow a prescribed schedule. A natural next step is to let the parameters depend on the observed trajectory, so that the extrapolation can react to the actual progress of the method. The analysis then has to handle parameter histories rather than a single sequence, and the dissipation coefficient of the energy is no longer automatically nonnegative. We show below that a signed energy, together with an explicit finite correction for its negative prefix, is sufficient.

Write the momentum coefficient as $\beta_k=(t_k-1)/t_{k+1}$, so that the extrapolated point is $y^k=x^k+\beta_k(x^k-x^{k-1})$. For the rational choice
\begin{equation}\label{eq:intro-rational}
 \beta_k=\frac{k-1}{k+\alpha-1},
 \qquad \alpha>0,
\end{equation}
the value $\alpha=3$, that is $\beta_k=(k-1)/(k+2)$ and $t_k=(k+1)/2$, is critical: the corresponding root parameters degenerate to $(a,b)=(0,0)$, and the rational rule is no longer an interior instance of the general theory. The three regimes behave differently:
\begin{itemize}[leftmargin=2.2em]
\item for $\alpha>3$ the root parameters satisfy $0<a<1$ and $a^2/4\le b\le1/4$, and the accelerated little-o rate is available;
\item for $\alpha=3$ the finite $O(N^{-2})$ bound survives, but the little-o rate is a separate question;
\item for $1<\alpha<3$ the root parameter $a=(\alpha-3)/(\alpha-1)$ is negative, the clock is still positive, and only a weaker polynomial bound follows from the signed energy.
\end{itemize}
Table~\ref{tab:rates} summarizes the rates established in this paper. It is the multiobjective composite counterpart of the scalar table in \cite{attouch2017subcritical}, with the additional feature that the momentum coefficient may depend on the whole history.

\begin{table}[!htbp]
\centering
\small
\setlength{\tabcolsep}{4pt}
\caption{Rates for the accelerated proximal gradient method with general stepsizes. The parameter $\alpha$ refers to the rational momentum coefficient \eqref{eq:intro-rational}; the root family is the history-dependent generalization of the interior range.}
\label{tab:rates}
\begin{tabular}{@{}P{0.22\textwidth}P{0.27\textwidth}P{0.23\textwidth}P{0.21\textwidth}@{}}
\toprule
Momentum rule & Merit bound & Displacement & Status\\
\midrule
root family, $\alpha>3$ &
$U_k=o(k^{-2})$ &
$\norm{v^k}=o(k^{-1})$ &
closed, any history\\
rational, $\alpha=3$ &
$U_N\le\dfrac{2\ell R}{N(N+3)}$;\newline
$U_k=o(k^{-2})\Leftrightarrow\Acrit=0$ &
$\norm{v^k}=O(N^{-1})$ &
critical energy\\
rational, $1<\alpha<3$ &
$U_N=O\bigl(N^{-(\alpha-1)}\bigr)$ &
not asserted &
coarse bound\\
rational, $\alpha>0$ &
$U_k\to0$ &
$\norm{v^k}\to0$ &
qualitative\\
\bottomrule
\end{tabular}
\end{table}

Two features of the vector problem require particular care and explain why the scalar arguments do not transfer verbatim. First, for a comparison point $z$, the signed comparison
\begin{equation}\label{eq:intro-sigma}
 \sigma_k(z)=\min_{1\le i\le m}\bigl\{F_i(x^k)-F_i(z)\bigr\}
\end{equation}
can be negative even when $z$ is weakly Pareto optimal. A scalar proof that discards terms containing a nonnegative objective gap therefore cannot be copied without an additional argument. Second, the merit function is a supremum over comparison points, so pointwise convergence of $k^2\sigma_k(z)$ for every fixed $z$ does not imply convergence of $k^2\sup_z\sigma_k(z)$. The estimates below keep signed comparisons until a bounded, iteration-independent family of weak Pareto value representatives is available, and only then take suprema.

The main contributions of this paper are the following. 

\begin{itemize}
    \item The first group of results concerns a two-parameter root momentum rule,
\begin{equation}\label{eq:intro-root}
 t_{k+1}=\sqrt{t_k^2-a_kt_k+b_k}+\tfrac12,
 \qquad
 \beta_k=\frac{t_k-1}{t_{k+1}},
\end{equation}
where $0\le a_k<1$, $a_k^2/4\le b_k\le1/4$, and $a_k\to\bar a\in(0,1)$. No convergence or bounded-variation condition is imposed on $b_k$. We prove boundedness, $\sum_k kD_k<\infty$, $U_k=o(k^{-2})$, and $\norm{v^k}=o(k^{-1})$ for every admissible history, where $D_k=\norm{x^k-x^{k-1}}^2$. The coefficient that governs weighted displacement summability need only become positive eventually; its finite negative prefix produces an explicit trajectory-dependent correction and does not restart the algorithm.

\item The second group concerns the critical parameter $\alpha=3$. For general convex composite data we establish the finite merit bound
\begin{equation}\label{eq:intro-critical}
 U_N\le\frac{2\ell R}{N(N+3)},
\end{equation}
square summability of the displacements, displacement decay $O(1/N)$, and convergence of the full sequence to a weak Pareto point. The argument identifies a nonnegative critical energy $\Acrit$ such that
\begin{equation}\label{eq:intro-equivalence}
 U_k=o(k^{-2})\quad\Longleftrightarrow\quad \Acrit=0.
\end{equation}
Thus the remaining little-o question at $\alpha=3$ is not another large-$O$ estimate but the exclusion of a positive limiting energy.

\item The third group concerns the noncritical ranges. For $1<\alpha<3$ we give the rate
\begin{equation}\label{eq:intro-subcritical}
 U_N=O\bigl(N^{-(\alpha-1)}\bigr),
\end{equation}
without claiming optimality. For every $\alpha>0$ we prove the qualitative conclusions $U_k\to0$, $D_k\to0$, $F_i(x^k)\to c_i\in\R$, and $\sum_k D_k/k<\infty$.

\item Finally, the proofs are organized so that the same small set of estimates serves every regime. The energy function, the coefficient asymptotics, and the uniform window estimate are stated once in \cref{sec:key-lemmas}; the regime-specific sections then only combine them. This modular arrangement is what allows the history-dependent rule and the rational rule to be treated by the same argument.
\end{itemize}

The organization of the paper is as follows.
\Cref{sec:prelim} fixes the notation, the ordering, the merit function, the standing geometric assumption, and the momentum rules. \Cref{sec:algorithm} introduces the offset subproblem, its optimality condition, and the two instances of the accelerated proximal gradient method with general stepsizes. \Cref{sec:key-lemmas} collects the key inequalities, the energy function, the coefficient asymptotics, and the auxiliary lemmas. \Cref{sec:alpha-gt3} proves the rate for the history-dependent root family, \cref{sec:alpha-eq3} treats the critical parameter $\alpha=3$, \cref{sec:alpha-between} gives the coarse rate for $1<\alpha<3$, and \cref{sec:alpha-positive} the qualitative result for every $\alpha>0$. \Cref{sec:numerical} illustrates the little-$o$ and $O(k^{-2})$ regimes numerically on the JOS1 problem, and \cref{sec:conclusion} concludes.

\section{Preliminaries}\label{sec:prelim}


We study $\min_{x\in\R^n}F(x)=\bigl(F_1(x),\ldots,F_m(x)\bigr),
 \qquad F_i=f_i+g_i,\qquad i=1,\ldots,m.$ \eqref{eq:intro-problem} with $m,n$ finite. Set
\begin{equation}\label{eq:X}
 X=\bigcap_{i=1}^m\dom g_i\ne\varnothing,
 \qquad
 L=\max_iL_i,
 \qquad
 \ell\ge L,\quad \ell>0,
\end{equation}
where $L_i$ is a global Lipschitz constant for $\nabla f_i$ and $1/\ell$ is the step size. Smooth terms remain evaluable outside $X$, whereas iterates produced by the proximal model lie in $X$.

For $u,v\in\R^m$ write $u\le v$ when $v-u\in\R^m_+$ and $u<v$ when $v-u\in\R^m_{++}$. A point $p\in X$ is weakly Pareto optimal if no $z\in X$ satisfies $F(z)<F(p)$. The weak Pareto set is denoted by $\Xstar$. The merit function, in the sense of \cite{tanabe2019proximal}, is
\begin{equation}\label{eq:merit}
 u_0(x)=\sup_{z\in X}\min_{1\le i\le m}\bigl\{F_i(x)-F_i(z)\bigr\},
 \qquad x\in X.
\end{equation}
Taking $z=x$ gives $u_0(x)\ge0$, and $x$ is weakly Pareto optimal exactly when $u_0(x)=0$. The value can be infinite at a general feasible point; finiteness along the algorithm is proved below. We abbreviate $U_k=u_0(x^k)$ and use the signed comparison
\begin{equation}\label{eq:sigma}
 \sigma_k(z)=\min_{1\le i\le m}\bigl\{F_i(x^k)-F_i(z)\bigr\},
 \qquad z\in X,
\end{equation}
so that $U_k=\sup_{z\in X}\sigma_k(z)$. Unlike $U_k$, the quantity $\sigma_k(z)$ may be negative; this is the main difference from the scalar case, and it is the reason why the estimates below keep $\sigma_k(z)$ signed until a bounded family of comparison points is available.

\begin{assumption}\label{ass:geometry}
Fix $x^0\in X$ and let $\mathcal L_0=\{x\in X:F(x)\le F(x^0)\}$ be the initial lower level set. Every $x\in\mathcal L_0$ is dominated by a weak Pareto point: there is $p\in\Xstar$ with $F(p)\le F(x)$. Moreover,
\begin{equation}\label{eq:radius}
 R:=\sup_{\xi\in F(\Xstar\cap\mathcal L_0)}
 \inf_{\substack{z\in X\\F(z)=\xi}}\norm{z-x^0}^2<\infty.
\end{equation}
\end{assumption}

Neither boundedness of $\mathcal L_0$ nor attainment of either extremum in \eqref{eq:radius} is assumed. The inner infimum may be approached rather than attained; representatives are chosen within an arbitrarily small error. The assumption controls representatives of weak Pareto objective values, not the whole initial lower level set. It is the standard geometric hypothesis in the convergence analysis of multiobjective proximal gradient methods \cite{tanabe2019proximal}.

With this geometric assumption in place, we next introduce the momentum coefficients that govern the extrapolation step and thus determine the iterates at which the above representative-selection condition is enforced. We study two parametric families of momentum coefficients. Both are written in the form $\beta_k=(t_k-1)/t_{k+1}$, and both extrapolate at the point
\begin{equation}\label{eq:extrapolation}
 y^k=x^k+\beta_k(x^k-x^{k-1}).
\end{equation}

\emph{History-dependent root momentum.} Choose $a_k,b_k$ from the already computed history and update
\begin{equation}\label{eq:root}
 t_{k+1}=\sqrt{t_k^2-a_kt_k+b_k}+\frac12,
 \qquad
 \beta_k=\frac{t_k-1}{t_{k+1}}.
\end{equation}
The pair $(a_k,b_k)$ is admissible when
\begin{equation}\label{eq:rootconditions}
 0\le a_k<1,
 \qquad
 \frac{a_k^2}{4}\le b_k\le\frac14,
 \qquad
 a_k\longrightarrow\bar a\in(0,1);
\end{equation}
these conditions are imposed in \cref{sec:alpha-gt3}, where the rate is proved. The second parameter $b_k$ need not converge and need not have bounded total variation, so \eqref{eq:root} is a genuine history-dependent rule and not a prescribed schedule.

\emph{Rational momentum.} For a fixed $\alpha>0$,
\begin{equation}\label{eq:rational}
 \beta_k=\frac{k-1}{k+\alpha-1},
 \qquad
 t_k=\frac{k+\alpha-2}{\alpha-1}.
\end{equation}
This is the specialization of \eqref{eq:root} with
\begin{equation}\label{eq:rational-map}
 a=\frac{\alpha-3}{\alpha-1},
 \qquad
 b=\frac{a^2}{4},
\end{equation}
so for $\alpha>3$ the pair $(a,b)$ lies in the interior of the admissible region \eqref{eq:rootconditions}, at $\alpha=3$ it degenerates to $(a,b)=(0,0)$, and for $1<\alpha<3$ it gives $a<0$. For $\alpha<1$ the displayed $t_k$ is negative and serves only as a formal scaling; the algorithm uses $\beta_k$ directly, and the qualitative result for every $\alpha>0$ in \cref{sec:alpha-positive} does not require $t_k\ge1$.

Throughout,
\begin{equation}\label{eq:notation}
 v^k=x^k-x^{k-1},
 \quad
 D_k=\norm{v^k}^2,
 \quad
 q_k=\frac{\ell}{2}D_k,
 \quad
 w_k=t_kx^k-(t_k-1)x^{k-1}.
\end{equation}
For an inertial clock define the coefficient sequence
\begin{equation}\label{eq:coeff}
\begin{aligned}
 \eta_k&=1-\beta_k^2,
 & \delta_k&=t_k^2-t_{k+1}^2+t_{k+1},\\
 S_k&=\sum_{j=1}^k\delta_j,
 & C_N&=\sum_{j=1}^N t_j=t_N^2+S_{N-1}.
\end{aligned}
\end{equation}
The signed dissipation coefficient of the energy is
\begin{equation}\label{eq:Pdef}
 P_k=\eta_kS_k-\delta_k
 =(1-\beta_k^2)\sum_{j=1}^{k+1}t_j-2t_k-t_{k+1}+1.
\end{equation}
Write $P_k^\pm=\max\{\pm P_k,0\}$ and
\begin{equation}\label{eq:Ncorr}
 \Ncorr_N=\sum_{k=1}^{N-1}P_k^-D_k,
 \qquad
 \Ncorr_\infty=\sum_{k\ge1}P_k^-D_k.
\end{equation}
The identity $C_N=t_N^2+S_{N-1}$ follows from $t_1=1$ and telescoping.

We use the convention that an empty sum is zero. Constants may depend on the fixed problem, the step size, and the particular parameter history; they are not claimed to be uniform over all histories.

\section{Accelerated proximal gradient methods with general stepsizes}\label{sec:algorithm}
In this section we introduce the accelerated proximal gradient method for multiobjective optimization with general stepsizes. The method is a dynamic extension of the multiobjective proximal gradient method of \cite{tanabe2022root}. The algorithm is described in \cref{sec:algorithm-body}, and its convergence properties are analyzed in \cref{sec:key-lemmas,sec:alpha-gt3}.

\subsection{Subproblem and optimality condition}\label{sec:subproblem}

For $x\in X$ and $y\in\R^n$, define
\begin{align}
 \phi_i(z;x,y)&=f_i(y)+\ip{\nabla f_i(y)}{z-y}+g_i(z)-F_i(x),\label{eq:phi}\\
 \psi(z;x,y)&=\max_{1\le i\le m}\phi_i(z;x,y),\nonumber\\
 p_\ell(x,y)&=\argmin_{z\in\R^n}\left\{\psi(z;x,y)+\frac{\ell}{2}\norm{z-y}^2\right\}.\label{eq:prox}
\end{align}
The offsets $-F_i(x)$ are retained separately inside the maximum. Removing them changes the model in general; \eqref{eq:prox} is not the proximal gradient step for $\max_iF_i$, and it is not obtained from a fixed scalarization.

\begin{lemma}\label{lem:wellposed}
For every $x\in X$ and $y\in\R^n$, the minimum in \eqref{eq:prox} is attained at a unique point in $X$. Writing $p=p_\ell(x,y)$,
\begin{equation}\label{eq:aggregate-opt}
 \ell(y-p)\in\partial\psi(p;x,y).
\end{equation}
\end{lemma}
\begin{proof}
The maximum model is proper, lower semicontinuous, and convex because it is finite at any point of $X$. A proper closed convex function on a finite-dimensional space has an affine minorant \cite{rockafellar1970convex}. Applying this to one component and adding $\ell\norm{z-y}^2/2$ shows that the objective in \eqref{eq:prox} is coercive. Hence the minimum is attained. The objective is $\ell$-strongly convex, so the minimizer is unique. The convex Fermat rule and the sum rule with the differentiable quadratic give \eqref{eq:aggregate-opt}.
\end{proof}

The aggregate inclusion requires no componentwise constraint qualification. If all $g_i$ are finite and continuously differentiable near $p$, it becomes
\begin{equation}\label{eq:active-opt}
\begin{aligned}
 \ell(y-p)&=\sum_i\lambda_i\bigl(\nabla f_i(y)+\nabla g_i(p)\bigr),\\
 \lambda&\in\Delta_m,
 \qquad
 \lambda_i>0\Longrightarrow\phi_i(p;x,y)=\psi(p;x,y),
\end{aligned}
\end{equation}
where $\Delta_m=\{\lambda\in\R^m:\lambda\ge0,\ \sum_i\lambda_i=1\}$. The weights are determined jointly with the minimizer, not prescribed in advance.

\subsection{Accelerated proximal gradient method}\label{sec:algorithm-body}

We initialize $x^{-1}=x^0$ and compute $x^1=p_\ell(x^0,x^0)$. For $k\ge1$, the inertial step is
\begin{equation}\label{eq:iteration}
 y^k=x^k+\beta_k(x^k-x^{k-1}),
 \qquad
 x^{k+1}=p_\ell(x^k,y^k),
 \qquad
 \beta_k=\frac{t_k-1}{t_{k+1}},
 \quad t_1=1.
\end{equation}
In particular, $\beta_1=0$: both the initial step and the next step are noninertial. All statements concern exact solutions of \eqref{eq:prox}; the word general refers to the step size $1/\ell$ with $\ell\ge L$, not to an inexact subproblem solver.

The momentum coefficient is the only quantity that distinguishes the two instances below. The first instance uses the root recurrence that generalizes the momentum factor of \cite{tanabe2022root} to a history-dependent choice; the second uses the rational coefficient whose $\alpha>3$ range was studied in \cite{zhang2023}. The two instances are summarized in \cref{alg:root} and \cref{alg:rational} and differ only in the update of $\beta_k$.

\begin{algorithm}[!htbp]
\renewcommand{\algorithmicrequire}{\textbf{Input:}}
\renewcommand{\algorithmicensure}{\textbf{Output:}}
\caption{Accelerated proximal gradient method with general stepsizes (history-dependent root momentum)}
\label{alg:root}
\begin{algorithmic}[1]
\REQUIRE $x^0\in X$, $\ell\ge L$, $\ell>0$, and an admissible parameter rule.
\STATE Set $x^{-1}=x^0$, $t_1=1$; compute $x^1=p_\ell(x^0,x^0)$.
\FOR{$k=1,2,\ldots$}
 \STATE Choose $a_k,b_k$ from the already computed history.
 \STATE Set $t_{k+1}=\sqrt{t_k^2-a_kt_k+b_k}+1/2$ and $\beta_k=(t_k-1)/t_{k+1}$.
 \STATE Set $y^k=x^k+\beta_k(x^k-x^{k-1})$.
 \STATE Solve \eqref{eq:prox} exactly to obtain $x^{k+1}=p_\ell(x^k,y^k)$.
\ENDFOR
\ENSURE The iterates $(x^k)$.
\end{algorithmic}
\end{algorithm}

\begin{algorithm}[!htbp]
\renewcommand{\algorithmicrequire}{\textbf{Input:}}
\renewcommand{\algorithmicensure}{\textbf{Output:}}
\caption{Accelerated proximal gradient method with general stepsizes (rational momentum)}
\label{alg:rational}
\begin{algorithmic}[1]
\REQUIRE $x^0\in X$, $\ell\ge L$, $\ell>0$, and $\alpha>0$.
\STATE Set $x^{-1}=x^0$, $t_1=1$; compute $x^1=p_\ell(x^0,x^0)$.
\FOR{$k=1,2,\ldots$}
 \STATE Set $\beta_k=(k-1)/(k+\alpha-1)$.
 \STATE Set $y^k=x^k+\beta_k(x^k-x^{k-1})$.
 \STATE Solve \eqref{eq:prox} exactly to obtain $x^{k+1}=p_\ell(x^k,y^k)$.
\ENDFOR
\ENSURE The iterates $(x^k)$.
\end{algorithmic}
\end{algorithm}

Without acceleration, $\beta_k\equiv0$ and the iteration reduces to the multiobjective proximal gradient method of \cite{tanabe2019proximal}. The root rule \eqref{eq:root} is the mechanism that produces $\beta_k\to1$ at the rate $1/k$ and thus the $o(k^{-2})$ merit decay; it is analysed in \cref{sec:key-lemmas} and \cref{sec:alpha-gt3}.


\section{Key inequalities and energy function}\label{sec:key-lemmas}

The estimates in this section depend only on the exact proximal model and the inertial identity $\beta_k=(t_k-1)/t_{k+1}$, except for the coefficient asymptotics, which is stated for the root clock and used in \cref{sec:alpha-gt3}. The remaining results are reused in every later convergence proof.


\begin{definition}[Energy function]\label{def:energy}
For $k\ge1$ and $z\in X$ define
\begin{equation}\label{eq:H}
 \HH_k(z)=
 \left(\sum_{j=1}^kt_j\right)\sigma_k(z)
 +\frac{\ell}{2}\norm{w_k-z}^2
 +\frac{\ell}{2}\left(\sum_{j=1}^kt_j-t_k^2\right)D_k,
\end{equation}
with $w_k=t_kx^k-(t_k-1)x^{k-1}$ and $\sigma_k(z)$ as in \eqref{eq:sigma}. The coefficient $P_k$ in \eqref{eq:Pdef} is the signed dissipation coefficient of $\HH_k$.
\end{definition}

The coefficient $P_k$ measures the signed amount by which the energy decreases at step $k$. It need not be nonnegative at every step. The corrected energy keeps the negative contributions and accumulates them in $\Ncorr_N$; this is the only role of the negative prefix. In the scalar case the analogous function is the standard Lyapunov function of the accelerated gradient method, and $\sigma_k(z)$ reduces to the objective gap.


The following theorem collects the estimates that are reused in every later section. It uses only the definition of the exact proximal model and the inertial identity $\beta_k=(t_k-1)/t_{k+1}$.

\begin{theorem}\label{thm:basic}
Let $x^{k+1}$ be generated by \eqref{eq:iteration} with $t_k\ge1$ and $\beta_k=(t_k-1)/t_{k+1}$. For $z\in X$ let $\sigma_k(z)=\min_i\{F_i(x^k)-F_i(z)\}$. Then for every $k\ge1$ and every $z\in X$:
\begin{enumerate}[label=(\roman*),leftmargin=3em]
\item strong convexity at the minimizer gives
\begin{equation}\label{eq:strong}
\begin{aligned}
 &\psi(x^{k+1};x^k,y^k)+\frac{\ell}{2}\norm{x^{k+1}-y^k}^2
 +\frac{\ell}{2}\norm{x^{k+1}-z}^2\\
 &\qquad\le\psi(z;x^k,y^k)+\frac{\ell}{2}\norm{z-y^k}^2;
\end{aligned}
\end{equation}
\item the model dominates the true componentwise increment:
\begin{equation}\label{eq:model-lower}
 \max_i\{F_i(x^{k+1})-F_i(x^k)\}
 \le\psi(x^{k+1};x^k,y^k)+\frac{\ell}{2}\norm{x^{k+1}-y^k}^2,
\end{equation}
and
\begin{equation}\label{eq:model-upper}
 \psi(z;x^k,y^k)\le-\sigma_k(z);
\end{equation}
\item the first basic inequality is
\begin{equation}\label{eq:component}
 F_i(x^{k+1})+\frac{\ell}{2}\norm{x^{k+1}-x^k}^2
 \le F_i(x^k)+\frac{\ell\beta_k^2}{2}\norm{x^k-x^{k-1}}^2,
\end{equation}
and consequently
\begin{equation}\label{eq:gap-monotone}
 \sigma_{k+1}(z)+\frac{\ell}{2}\norm{x^{k+1}-x^k}^2
 \le\sigma_k(z)+\frac{\ell\beta_k^2}{2}\norm{x^k-x^{k-1}}^2;
\end{equation}
equivalently, with $q_k=\ell D_k/2$ and $\eta_k=1-\beta_k^2$,
\begin{equation}\label{eq:augmented}
 \sigma_k(z)+q_k
 \ge\sigma_{k+1}(z)+q_{k+1}+\eta_kq_k;
\end{equation}
\item the second basic inequality is
\begin{equation}\label{eq:distance}
 \sigma_{k+1}(z)
 \le\frac{\ell}{2}\left(\norm{y^k-z}^2-\norm{x^{k+1}-z}^2\right).
\end{equation}
\end{enumerate}
\end{theorem}
\begin{proof}
(i) The objective in \eqref{eq:prox} is $\ell$-strongly convex and $x^{k+1}$ is its minimizer. The strong-convexity subgradient inequality gives \eqref{eq:strong}.

(ii) The descent lemma applied to each $f_i$ gives
\[
 f_i(x^{k+1})\le f_i(y^k)+\ip{\nabla f_i(y^k)}{x^{k+1}-y^k}
 +\frac{L_i}{2}\norm{x^{k+1}-y^k}^2.
\]
Adding $g_i(x^{k+1})-F_i(x^k)$, using $L_i\le\ell$, and taking the maximum over $i$ yields \eqref{eq:model-lower}. Convexity of $f_i$ gives
\[
 f_i(y^k)+\ip{\nabla f_i(y^k)}{z-y^k}\le f_i(z),
\]
and therefore $\psi(z;x^k,y^k)\le\max_i\{F_i(z)-F_i(x^k)\}=-\sigma_k(z)$, which is \eqref{eq:model-upper}. Taking $z=x^k$ gives $\psi(x^k;x^k,y^k)\le0$.

(iii) Take $z=x^k$ in \eqref{eq:strong} and replace the model terms using \eqref{eq:model-lower} and \eqref{eq:model-upper}. Since $y^k-x^k=\beta_k(x^k-x^{k-1})$,
\[
 \max_i\{F_i(x^{k+1})-F_i(x^k)\}
 +\frac{\ell}{2}\norm{x^{k+1}-x^k}^2
 \le\frac{\ell\beta_k^2}{2}\norm{x^k-x^{k-1}}^2.
\]
Subtracting the same $F_i(z)$ from both sides, taking the minimum over $i$, and noting that the displacement term does not depend on $i$ gives \eqref{eq:gap-monotone}.

(iv) For real sequences $(A_i)$ and $(B_i)$ one has
$\min_iA_i-\min_iB_i\le\max_i(A_i-B_i)$.
Taking $A_i=F_i(x^{k+1})-F_i(z)$ and $B_i=F_i(x^k)-F_i(z)$ gives
$\sigma_{k+1}(z)-\sigma_k(z)\le\max_i\{F_i(x^{k+1})-F_i(x^k)\}$.
Combining \eqref{eq:strong}, \eqref{eq:model-lower}, and \eqref{eq:model-upper} yields
\[
 \max_i\{F_i(x^{k+1})-F_i(x^k)\}
 \le-\sigma_k(z)+\frac{\ell}{2}\left(\norm{y^k-z}^2-\norm{x^{k+1}-z}^2\right).
\]
The two inequalities cancel $-\sigma_k(z)$ and give \eqref{eq:distance}.
\end{proof}

Although we cannot guarantee that the function values at the iterates decrease monotonically, we can ensure that they never exceed the initial value. Since $\beta_1=0$ and $y^1=x^1$, applying \cref{thm:basic}(iii)--(iv) with $x^0$ in place of both $x^k$ and $y^k$ gives, for every $z\in X$,
\begin{equation}\label{eq:initial}
\begin{aligned}
 F_i(x^1)+\frac{\ell}{2}\norm{x^1-x^0}^2&\le F_i(x^0),\\
 \sigma_1(z)+\frac{\ell}{2}\norm{x^1-z}^2&\le\frac{\ell}{2}\norm{x^0-z}^2.
\end{aligned}
\end{equation}
Then \eqref{eq:component} and $\beta_k^2\le1$ give
\[
 F_i(x^{k+1})+\frac{\ell}{2}D_{k+1}
 \le F_i(x^k)+\frac{\ell}{2}D_k
 \le\cdots\le F_i(x^0),
\]
and dropping the nonnegative square term yields the invariance of the initial lower level set,
\begin{equation}\label{eq:level}
 F(x^k)\le F(x^0),
 \qquad
 x^k\in\mathcal L_0,
 \qquad k\ge0.
\end{equation}

Next, we introduce the energy function and prove its basic properties. 

\begin{lemma}\label{lem:energy}
Let $W_k=\sum_{j=1}^kt_j-t_k^2$. Then $W_k\ge0$, $W_1=0$, and for every $k,N\ge1$ and $z\in X$,
\begin{equation}\label{eq:energy-recursion}
 \HH_{k+1}(z)\le\HH_k(z)-\frac{\ell}{2}P_kD_k,
\end{equation}
and
\begin{equation}\label{eq:cumulative}
\begin{aligned}
 \HH_N(z)+\frac{\ell}{2}\sum_{k=1}^{N-1}P_k^+D_k
 \le\frac{\ell}{2}\left(\norm{x^0-z}^2+\sum_{k=1}^{N-1}P_k^-D_k\right).
\end{aligned}
\end{equation}
\end{lemma}
\begin{proof}
Squaring the root recurrence after moving the square root gives
\[
 t_k^2-t_{k+1}^2+t_{k+1}=a_kt_k+\frac14-b_k\ge0.
\]
Summing from $j=1$ to $k-1$ and using $t_1=t_1^2=1$ gives
$W_k=\sum_{j=1}^{k-1}(a_jt_j+1/4-b_j)\ge0$ and $W_1=0$.

Multiply \eqref{eq:gap-monotone} by $t_{k+1}-1$, add \eqref{eq:distance}, and multiply by $t_{k+1}$. Using
\[
 t_{k+1}y^k-(t_{k+1}-1)x^k=t_kx^k-(t_k-1)x^{k-1}
\]
and expanding squared norms gives
\begin{equation}\label{eq:accelerated}
\begin{aligned}
 &t_{k+1}^2\sigma_{k+1}(z)
 +\frac{\ell}{2}\norm{t_{k+1}x^{k+1}-(t_{k+1}-1)x^k-z}^2\\
 &\quad\le t_{k+1}(t_{k+1}-1)\sigma_k(z)
 +\frac{\ell}{2}\norm{t_kx^k-(t_k-1)x^{k-1}-z}^2.
\end{aligned}
\end{equation}
Add $W_{k+1}\ge0$ times \eqref{eq:gap-monotone}. The comparison coefficients become
$t_{k+1}^2+W_{k+1}=\sum_{j=1}^{k+1}t_j$ and
$t_{k+1}(t_{k+1}-1)+W_{k+1}=\sum_{j=1}^kt_j$, and the displacement coefficient difference is
\[
 W_k-\beta_k^2W_{k+1}=P_k.
\]
This proves \eqref{eq:energy-recursion}. Summing from $k=1$ to $N-1$, using \eqref{eq:initial}, and splitting $P_k=P_k^+-P_k^-$ gives \eqref{eq:cumulative}.
\end{proof}

The next proposition gives the asymptotic behavior of the coefficients and weights for the root clock. It is used in \cref{sec:alpha-gt3} to prove the rate. The proof is a straightforward application of the recurrence \eqref{eq:root} and the admissibility conditions \eqref{eq:rootconditions}.

\begin{proposition}\label{prop:limits}
Assume \eqref{eq:rootconditions} and put $c=(1-\bar a)/2\in(0,1/2)$. Then
\begin{enumerate}[label=(\roman*),leftmargin=3em]
\item $t_{k+1}-t_k\to c$ and $t_k/k\to c$;
\item $\beta_k\to1$ and $k(1-\beta_k^2)\to2(1+c)/c$;
\item $S_k/k^2\to\bar ac/2$ and $C_N/N^2\to c/2$;
\item $P_k/k\to\bar a$. In particular, there is $k_0\ge2$ such that
\begin{equation}\label{eq:Ppositive}
 P_k\ge\frac{\bar a}{2}k>0,
 \qquad k\ge k_0.
\end{equation}
\end{enumerate}
\end{proposition}
\begin{proof}
Rationalize the recurrence:
\[
 t_{k+1}-t_k=\frac{1-a_k}{2}
 +\frac{b_k-a_k^2/4}{\sqrt{(t_k-a_k/2)^2+b_k-a_k^2/4}+t_k-a_k/2}.
\]
The fraction is bounded by $1/[4(2t_k-1)]$. Since $a_k\to\bar a<1$, eventually $1-a_k\ge(1-\bar a)/2$, so $t_{k+1}-t_k\ge(1-\bar a)/4>0$. Thus $t_k\to\infty$, the fraction tends to zero, and $t_{k+1}-t_k\to c$. Ces\`aro averaging gives $t_k/k\to c$.

Since $\beta_k=(t_k-1)/t_{k+1}$,
\[
 k(1-\beta_k)=\frac{k(1+t_{k+1}-t_k)}{t_{k+1}}
 =\frac{1+t_{k+1}-t_k}{t_{k+1}/k}
 \longrightarrow\frac{1+c}{c}.
\]
Hence $\beta_k\to1$ and $k(1-\beta_k^2)=k(1-\beta_k)(1+\beta_k)\to2(1+c)/c$.

For (iii), fix $0<\delta<c$. For large $j$, $(c-\delta)j\le t_j\le(c+\delta)j$. Dividing the partial sums by $k^2$ and using $\sum_{j=J}^{k}j/k^2\to1/2$ gives $S_k/k^2\to\bar ac/2$ because $\delta_k/k\to\bar ac$; the same argument with $S_k=\sum_{j=1}^{k+1}t_j-t_{k+1}^2$ gives $C_N/N^2\to c/2$.

Finally,
\[
 \frac{P_k}{k}
 =(k\eta_k)\frac{S_k}{k^2}-\frac{\delta_k}{k}
 \longrightarrow\frac{2(1+c)}{c}\cdot\frac{\bar ac}{2}-\bar ac
 =1-2c=\bar a>0.
\]
Thus \eqref{eq:Ppositive} holds eventually. Negative $P_k$ can occur only in a finite prefix.
\end{proof}


The four lemmas below depend only on the basic estimates, the initial lower level set, and the coefficient asymptotics. They are proved before the main argument so that the final proof contains no repeated algebra. \cref{lem:ballsup} restricts the comparison points to a ball of radius $R+\varepsilon$ around the initial point. 

\begin{lemma}\label{lem:ballsup}
For every $N\ge1$ and every $\varepsilon>0$,
\begin{equation}\label{eq:ballsup}
 u_0(x^N)=\sup_{\substack{z\in X\\\norm{z-x^0}^2\le R+\varepsilon}}\sigma_N(z).
\end{equation}
\end{lemma}
\begin{proof}
Let $z\in X$ satisfy $F(z)\le F(x^N)$. By \eqref{eq:level}, $z\in\mathcal L_0$. By \cref{ass:geometry} there is $p\in\Xstar$ with $F(p)\le F(z)\le F(x^N)\le F(x^0)$, hence $p\in\Xstar\cap\mathcal L_0$. The definition of $R$ gives, for every $\varepsilon>0$, a representative $\widehat z$ with $F(\widehat z)=F(p)$ and $\norm{\widehat z-x^0}^2\le R+\varepsilon$. Since $F_i(\widehat z)\le F_i(z)$ for every $i$, we have $\sigma_N(\widehat z)\ge\sigma_N(z)$. Thus every comparison is dominated by one in the ball. Applying this to $z=x^N$ gives a point in the ball with nonnegative comparison, so the restricted supremum is at least zero. Comparisons with a negative value do not affect the supremum, and the reverse inequality is immediate.
\end{proof}

\cref{lem:bounded} gives a uniform bound on the orbit. 

\begin{lemma}\label{lem:bounded}
Let $k_0$ be given by \eqref{eq:Ppositive} and set
\[
 B=\left[4(R+1)+2\sum_{k=1}^{k_0-1}P_k^-D_k\right]^{1/2}<\infty.
\]
Then
\begin{equation}\label{eq:orbitbound}
 \norm{x^N-x^0}\le B,
 \qquad N\ge0.
\end{equation}
\end{lemma}
\begin{proof}
In \eqref{eq:cumulative}, keep the squared-distance term of the combined point, drop the nonnegative positive-coefficient sum and the nonnegative $D_N$ term, and replace the comparison point center by $x^0$ using
$\norm{w_N-x^0}^2\le2\norm{w_N-z}^2+2\norm{z-x^0}^2$. Taking the supremum over $\norm{z-x^0}^2\le R+1$ and using \cref{lem:ballsup} gives
\[
 \norm{w_N-x^0}^2\le4(R+1)+2\sum_{k=1}^{k_0-1}P_k^-D_k.
\]
The right side is independent of $N$. Since
\[
 x^N-x^0=\frac1{t_N}(w_N-x^0)+\left(1-\frac1{t_N}\right)(x^{N-1}-x^0)
\]
is a convex combination because $t_N\ge1$, induction gives \eqref{eq:orbitbound}.
\end{proof}

\cref{lem:monotone} shows that the merit function is nonincreasing.

\begin{lemma}\label{lem:monotone}
For every $k\ge1$,
\begin{equation}\label{eq:meritmonotone}
\begin{aligned}
 0\le{}&U_{k+1}+\frac{\ell}{2}D_{k+1}\\
 \le{}&U_k+\frac{\ell\beta_k^2}{2}D_k
 \le U_k+\frac{\ell}{2}D_k.
\end{aligned}
\end{equation}
\end{lemma}
\begin{proof}
Inequality \eqref{eq:gap-monotone} holds for every comparison point. Taking the supremum, using $\beta_k^2\le1$, and noting that the displacement term is independent of $z$ gives \eqref{eq:meritmonotone}.
\end{proof}

\cref{lem:tail} gives a uniform estimate that is used in every later section.

\begin{lemma}\label{lem:tail}
For every integer $K\ge1$ there are
$r=r(K)\in\{K,\ldots,2K\}$ and $s=s(K)\in\{3K,\ldots,4K\}$ such that
\begin{equation}\label{eq:tailbound}
\begin{aligned}
 0\le{}&K^2\left(U_s+\frac{\ell}{2}D_s\right)\\
 \le{}&\frac{\ell}{2}\left[
 C(\sqrt2+2+\sqrt{\log4})\left(\sum_{j=K}^{4K}jD_j\right)^{1/2}
 +3\sum_{j=K}^{4K}jD_j\right],
\end{aligned}
\end{equation}
where
\begin{equation}\label{eq:constant}
 C=2(B+\sqrt{R+1})\left[1+\sup_{\nu\ge1}\nu(1-\beta_\nu)\right]<\infty.
\end{equation}
\end{lemma}
\begin{proof}
For $z\in X$ with $\norm{z-x^0}^2\le R+1$, \eqref{eq:orbitbound} gives $\norm{x^j-z}\le B+\sqrt{R+1}$. Let $h_j(z)=\norm{x^j-z}^2$ and $\Delta h_j=h_j-h_{j-1}$. Expanding \eqref{eq:distance} gives
\begin{equation}\label{eq:distance-expansion}
 \frac2\ell\sigma_{k+1}(z)
 \le\beta_k\Delta h_k-\Delta h_{k+1}+(\beta_k+\beta_k^2)D_k.
\end{equation}
Choose $r\in\{K,\ldots,2K\}$ minimizing $jD_j$ and $s\in\{3K,\ldots,4K\}$ minimizing $jD_j$ on the respective intervals. Then
\[
 r^2D_r\le2\sum_{j=K}^{4K}jD_j,
 \qquad
 s^2D_s\le4\sum_{j=K}^{4K}jD_j,
\]
and these indices do not depend on $z$. Add $D_{k+1}$ to \eqref{eq:distance-expansion}, use monotonicity of $\sigma_j(z)+q_j$, multiply by $k$, and sum from $r$ to $s-1$. Summation by parts gives
\begin{align}
 \frac2\ell\left(\sum_{k=r}^{s-1}k\right)(\sigma_s(z)+q_s)
 &\le r\beta_r\Delta h_r-(s-1)\Delta h_s\nonumber\\
 &\quad+\sum_{j=r+1}^{s-1}[j\beta_j-(j-1)]\Delta h_j\label{eq:windowsum}\\
 &\quad+\sum_{k=r}^{s-1}k\bigl[(\beta_k+\beta_k^2)D_k+D_{k+1}\bigr].\nonumber
\end{align}
Since $|\Delta h_j|\le2(B+\sqrt{R+1})\sqrt{D_j}$ and $|j\beta_j-(j-1)|\le1+\sup_\nu\nu(1-\beta_\nu)$, the endpoint and interior distance terms are bounded by
$C(\sqrt2+2)\sqrt{T_K}$ and $C\sqrt{T_K\log4}$, respectively, with $T_K=\sum_{j=K}^{4K}jD_j$. The last sum is at most $3T_K$. Thus
\[
 \frac2\ell\left(\sum_{k=r}^{s-1}k\right)(\sigma_s(z)+q_s)
 \le C(\sqrt2+2+\sqrt{\log4})\sqrt{T_K}+3T_K.
\]
The right side is independent of $z$. Taking the supremum over the ball and using \cref{lem:ballsup} replaces $\sigma_s(z)$ by $U_s$. Since $\sum_{k=r}^{s-1}k\ge K^2$, the left side is at least $2K^2(U_s+q_s)/\ell$. Multiplying by $\ell/2$ gives \eqref{eq:tailbound}.
\end{proof}

\section{Convergence rate for \texorpdfstring{$\alpha>3$}{alpha>3}}\label{sec:alpha-gt3}

This section analyses \cref{alg:root} under the history-dependent conditions \eqref{eq:rootconditions}. The second parameter $b_k$ need not converge and need not have bounded total variation, and every conclusion holds for each admissible history separately. The modular chain of \cref{sec:key-lemmas} supplies the basic estimates, the initial level set, the corrected energy, the coefficient asymptotics, and the auxiliary lemmas; the final argument needs only to combine them. The rational rule with $\alpha>3$ is the special case \eqref{eq:rational-map} of the same construction, so the root rule is the history-dependent extension of the interior regime.

\begin{theorem}[History-dependent root momentum]\label{thm:root}
Under \cref{ass:geometry}, suppose \cref{alg:root} uses parameters satisfying \eqref{eq:rootconditions}. Then $(x^k)$ is bounded and
\begin{equation}\label{eq:rootrates}
 \sum_{k\ge1}kD_k<\infty,
 \qquad
 U_k=o(k^{-2}),
 \qquad
 \norm{v^k}=o(k^{-1}).
\end{equation}
Moreover, $P_k^-$ has finite support, $\Ncorr_\infty<\infty$, and
\begin{equation}\label{eq:rootfinite}
 0\le U_N\le\frac{\ell(R+\Ncorr_N)}{2C_N},
 \qquad
 \frac{C_N}{N^2}\longrightarrow\frac{1-\bar a}{4}.
\end{equation}
No convergence or bounded-variation assumption on $b_k$ is needed. The negative-prefix budget and its final index may depend on the history.
\end{theorem}

\begin{proof}
Take $\varepsilon>0$ and restrict to $\norm{z-x^0}^2\le R+\varepsilon$ in \eqref{eq:cumulative}. Dropping the nonnegative term $(\sum_{j=1}^{N-1}P_k^+D_k)$ and the nonnegative squared-distance term, and using \cref{lem:ballsup}, gives
\[
 C_NU_N
 \le\frac{\ell}{2}\left(R+\varepsilon+\Ncorr_N\right).
\]
Letting $\varepsilon\downarrow0$ proves the finite bound in \eqref{eq:rootfinite}. Keeping the positive-coefficient sum and deleting $C_NU_N$ gives the budget
\begin{equation}\label{eq:budget}
 \sum_{k=1}^{N-1}P_k^+D_k\le R+\Ncorr_N.
\end{equation}
For $N>k_0$, \eqref{eq:Ppositive} and \eqref{eq:budget} imply
\[
 \frac{\bar a}{2}\sum_{k=k_0}^{N-1}kD_k
 \le\sum_{k=k_0}^{N-1}P_kD_k
 \le R+\sum_{k=1}^{k_0-1}P_k^-D_k.
\]
The right side is independent of $N$, so $\sum_k kD_k<\infty$.

By \eqref{eq:meritmonotone}, $W_k=U_k+q_k$ is nonincreasing and nonnegative. Let $T_K=\sum_{j=K}^{\infty}jD_j\to0$. Then $0\le\sum_{j=K}^{4K}jD_j\le T_K\to0$. \Cref{lem:tail} gives $K^2W_{s(K)}\to0$. For $N\ge4$, take $K=\lfloor N/4\rfloor$; then $s(K)\le N$, $N/K\le7$, and monotonicity gives
\[
 0\le N^2W_N\le49K^2W_{s(K)}\longrightarrow0.
\]
Thus $N^2U_N\to0$ and $N^2D_N\to0$, which is \eqref{eq:rootrates}. Boundedness is \cref{lem:bounded}. The limit $C_N/N^2\to(1-\bar a)/4$ is \cref{prop:limits}(iii) with $c=(1-\bar a)/2$.
\end{proof}

\begin{remark}[What the proof uses]\label{rem:criterion}
The argument applies to any increasing clock $t_k\ge1$ with $t_1=1$, $\delta_k\ge0$, $\sup_k k(1-\beta_k)<\infty$, and $P_k\ge\gamma k$ eventually for some $\gamma>0$. The root recurrence is one sufficient construction. These algebraic conditions isolate the mechanism: an eventually positive dissipation and a momentum coefficient that approaches one at the rate $1/k$.
\end{remark}

\Cref{thm:root} covers the rational rule with $\alpha>3$ by taking $a_k\equiv a=(\alpha-3)/(\alpha-1)$ and $b_k\equiv a^2/4$, in which case $\Ncorr_N=0$. The endpoint $\alpha=3$ is excluded because it corresponds to $(a,b)=(0,0)$, which lies on the boundary of \eqref{eq:rootconditions}; it is analysed separately in \cref{sec:alpha-eq3}.

\section{Convergence rate for \texorpdfstring{$\alpha=3$}{alpha=3}}\label{sec:alpha-eq3}

In this section we set
\begin{equation}\label{eq:criticalclock}
\begin{aligned}
 \beta_k&=\frac{k-1}{k+2},
 & t_k&=\frac{k+1}{2},\\
 y^k&=x^k+\beta_k(x^k-x^{k-1}),
 & x^{k+1}&=p_\ell(x^k,y^k).
\end{aligned}
\end{equation}
This is the rational rule \eqref{eq:rational} with $\alpha=3$. It is not an instance of the interior condition \eqref{eq:rootconditions}: the corresponding root parameters are $(a,b)=(0,0)$, which lie on the boundary of the admissible region. The analysis therefore uses the model inequalities and the algebraic energy estimate directly. We keep the modular order used above: main results, energy and coefficients, proof of the finite bound, point convergence, and the critical energy.

\subsection{Algorithm and main results}

The iteration is \cref{alg:rational} with $\alpha=3$. The first result is the critical baseline.

\begin{theorem}\label{thm:critical}
Under \cref{ass:geometry}, the sequences $(x^k)$ and $(w_k)$ generated by \eqref{eq:criticalclock} are bounded, and
\begin{equation}\label{eq:criticalbaseline}
 U_N\le\frac{2\ell R}{N(N+3)},
 \qquad
 \sum_{k\ge1}D_k<\infty,
 \qquad
 \norm{v^k}=O(k^{-1}).
\end{equation}
\end{theorem}

The second result describes the limit and the critical energy.

\begin{theorem}\label{thm:criticalenergy}
Under \cref{ass:geometry}, there are finite numbers $c_i=\lim_kF_i(x^k)$ such that the set
\begin{equation}\label{eq:limitset}
 \mathcal C=\{z\in X:F_i(z)\le c_i\ \text{for every }i\}
\end{equation}
is nonempty, closed, and convex. The iterates converge to a point $\bar x\in\mathcal C\cap\Xstar$. There is a constant $\Acrit\ge0$ such that
\begin{equation}\label{eq:criticalmerit}
 t_k^2\bigl(\sigma_k(\bar x)+q_k\bigr)\longrightarrow\Acrit,
 \qquad
 t_k^2(U_k+q_k)\longrightarrow\Acrit,
 \qquad
 U_k-\sigma_k(\bar x)=o(k^{-2}).
\end{equation}
Moreover,
\begin{equation}\label{eq:criticalequivalence}
 U_k=o(k^{-2})
 \quad\Longleftrightarrow\quad
 \Acrit=0.
\end{equation}
\end{theorem}

The proof of \cref{thm:critical} occupies \cref{sec:critical-energy}--\cref{sec:critical-proof}; \cref{thm:criticalenergy} is proved in \cref{sec:point}--\cref{sec:critical-identification}. The equivalence in \eqref{eq:criticalequivalence} is the precise sense in which the critical little-o question is reduced to the vanishing of a nonnegative energy constant.

\subsection{Energy identity and coefficient asymptotics}\label{sec:critical-energy}

The basic estimates of \cref{thm:basic} use only the model definition and the identity $\beta_k=(t_k-1)/t_{k+1}$. They therefore apply verbatim to the critical clock. We record the specialized coefficients:
\begin{equation}\label{eq:criticalcoeff}
 \beta_k=\frac{k-1}{k+2},
 \qquad
 \eta_k=1-\beta_k^2=\frac{6k+3}{(k+2)^2},
 \qquad
 \delta_k=\frac14.
\end{equation}
In the critical case $q_k=\ell D_k/2$, and the componentwise inequality \eqref{eq:component} reads
\begin{equation}\label{eq:criticalcomponent}
 F_i(x^{k+1})+q_{k+1}\le F_i(x^k)+\frac{\ell\beta_k^2}{2}D_k.
\end{equation}

\begin{lemma}\label{lem:criticalcoeff}
For the critical clock,
\begin{equation}\label{eq:criticalP}
 \delta_k=\frac14,
 \qquad
 S_k=\frac{k}{4},
 \qquad
 C_N=\frac{N(N+3)}{4},
 \qquad
 P_k=\frac{(k-1)(5k+4)}{4(k+2)^2}.
\end{equation}
In particular, $P_1=0$ and $P_k\ge7/32$ for $k\ge2$.
\end{lemma}
\begin{proof}
The formulas for $\delta_k$ and $S_k$ are immediate. Since $C_N=t_N^2+S_{N-1}$,
\[
 C_N=\frac{(N+1)^2}{4}+\frac{N-1}{4}=\frac{N^2+3N}{4}.
\]
For $P_k$,
\[
 P_k=\eta_kS_k-\delta_k
 =\frac{6k+3}{(k+2)^2}\cdot\frac{k}{4}-\frac14
 =\frac{(k-1)(5k+4)}{4(k+2)^2}.
\]
Finally,
\[
 32P_k-7=\frac{3(k-2)(11k+10)}{(k+2)^2}\ge0,
 \qquad k\ge2,
\]
and $P_1=0$.
\end{proof}

Since $P_k\ge0$ for every $k$, the energy function of \cref{def:energy} satisfies the clean inequality
\begin{equation}\label{eq:criticalenergybound}
 \HH_N(z)\le\frac{\ell}{2}\norm{x^0-z}^2,
 \qquad
 \HH_N(z)=C_N\sigma_N(z)+\frac{\ell}{2}\norm{w_N-z}^2
 +\frac{\ell}{2}\left(C_N-t_N^2\right)D_N.
\end{equation}
This is \eqref{eq:cumulative} with $\Ncorr_N=0$.

\subsection{Proof of the critical baseline}\label{sec:critical-proof}

\begin{proof}[Proof of \Cref{thm:critical}]
In \eqref{eq:criticalenergybound}, keep the term $C_N\sigma_N(z)$, drop the two nonnegative squared-distance terms, restrict to $\norm{z-x^0}^2\le R+\varepsilon$, and take the supremum using \cref{lem:ballsup}. This gives
\[
 C_NU_N\le\frac{\ell}{2}(R+\varepsilon).
\]
Letting $\varepsilon\downarrow0$ and using $C_N=N(N+3)/4$ proves the merit bound. For the displacement sum, choose a representative $z_N$ with $F(z_N)\le F(x^N)$ and $\norm{z_N-x^0}^2\le R+1$. Then $\sigma_N(z_N)\ge0$, and keeping the displacement term in \eqref{eq:criticalenergybound} gives
\[
 \frac{\ell}{2}\left(C_N-t_N^2\right)D_N\le\frac{\ell}{2}\norm{x^0-z_N}^2
 \le\frac{\ell}{2}(R+1).
\]
Since $C_N-t_N^2=(N-1)/4$, we obtain $D_N=O(1/N)$ and hence $\sum_kD_k<\infty$.

Boundedness of $w_N$ uses the same representative: the first and third terms in \eqref{eq:criticalenergybound} are nonnegative, so
\[
 \norm{w_N-z_N}^2\le\norm{x^0-z_N}^2\le R+1,
\]
and therefore $\norm{w_N-x^0}\le2\sqrt{R+1}$. Since $x^N=t_N^{-1}w_N+(1-t_N^{-1})x^{N-1}$ and $t_N\ge1$, induction gives $\norm{x^N-x^0}\le2\sqrt{R+1}$. Finally,
\[
 w_k-x^k=(t_k-1)v^k,
\]
so
\[
 \norm{v^k}\le\frac{\norm{w_k}+\norm{x^k}}{t_k-1}=O(k^{-1}).
\]
\end{proof}

\subsection{Point convergence}\label{sec:point}

\begin{lemma}\label{lem:criticalpoint}
There are finite numbers $c_i=\lim_kF_i(x^k)$ such that $\mathcal C$ in \eqref{eq:limitset} is nonempty, closed, and convex. The iterates converge to a point $\bar x\in\mathcal C\cap\Xstar$. For every $z\in\mathcal C$, the energy $\EE_k(z)=t_k^2\sigma_k(z)+\frac{\ell}{2}\norm{w_k-z}^2$ converges and $\sum_k|\sigma_k(z)|<\infty$. With $\mathcal V=\operatorname{span}(\mathcal C-\mathcal C)$,
\begin{equation}\label{eq:projection}
 \ip{w_k-\bar x}{d}\longrightarrow0,
 \qquad d\in\mathcal V.
\end{equation}
\end{lemma}
\begin{proof}
By \eqref{eq:component} with $\beta_k^2\le1$, $F_i(x^k)+q_k$ is nonincreasing. It is bounded below because the orbit is bounded and each $F_i$ is proper and lower semicontinuous. Since $q_k\to0$, it has a finite limit $c_i$, and $F_i(x^k)\to c_i$. Lower semicontinuity places every cluster point in $\mathcal C$, so $\mathcal C\ne\varnothing$. It is closed and convex as an intersection of closed convex sublevel sets.

For $z\in\mathcal C$, monotonicity gives $\sigma_k(z)+q_k\ge0$. The energies are bounded below because
\[
 \EE_k(z)\ge-t_k^2q_k\ge-C
\]
by \eqref{eq:criticalbaseline}. Let $\epsilon_K=\frac14\sum_{j=K}^\infty q_j\to0$. From \eqref{eq:energy-recursion} with $\delta_k=1/4$,
\[
 \EE_{k+1}(z)\le\EE_k(z)-\frac14(\sigma_k(z)+q_k)+\frac14q_k.
\]
Hence $\EE_k(z)+\epsilon_k$ is nonincreasing and bounded below, so it converges and $\sum_k(\sigma_k(z)+q_k)<\infty$. Since $\sum_kq_k<\infty$, also $\sum_k|\sigma_k(z)|<\infty$.

For $z\in\mathcal C$ let $I(z)=\{i:F_i(z)=c_i\}$. This set is nonempty: otherwise $\lim_k\sigma_k(z)=\min_i(c_i-F_i(z))>0$, contradicting $\sigma_k(z)\to0$. If a convex $F_i\le c_i$ on $\mathcal C$ attains $c_i$ at a relative interior point, convexity along line segments through that point forces $F_i\equiv c_i$ on $\mathcal C$. Therefore
\[
 I(z)=I_0:=\{i:F_i\equiv c_i\text{ on }\mathcal C\}\ne\varnothing
 \qquad(z\in\ri\mathcal C).
\]
For fixed $z$, the positive gaps $c_i-F_i(z)$ for $i\notin I(z)$ eventually exclude those indices from the minimum, and
\begin{equation}\label{eq:activeC}
 \sigma_k(z)=\min_{i\in I(z)}\bigl(F_i(x^k)-c_i\bigr)
\end{equation}
for all sufficiently large $k$. For $z,z'\in\ri\mathcal C$, the comparison sequences eventually agree because both reduce to the same minimum over $I_0$. Convergence of $\EE_k(z)-\EE_k(z')$ then implies convergence of $\ip{w_k}{z-z'}$.

Differences of relative interior points span $\mathcal V$. Choose a finite basis among these differences; the projection of $w_k$ onto $\mathcal V$ converges. The identity
\[
 x^k=t_k^{-1}w_k+(1-t_k^{-1})x^{k-1},
 \qquad
 \sum_kt_k^{-1}=\infty,
\]
and the elementary averaging lemma imply that the projection of $x^k$ has the same limit. Any two cluster points lie in $\mathcal C$, and their difference belongs to $\mathcal V$; their projections coincide, so they are equal. Since the orbit is bounded, $x^k\to\bar x\in\mathcal C$ and \eqref{eq:projection} holds. If a feasible point strictly improved every component at $\bar x$, then $\min_i(c_i-F_i(z))>0$, contradicting $U_k\to0$. Hence $\bar x\in\Xstar$.
\end{proof}

\subsection{The critical energy}\label{sec:critical-identification}

\begin{lemma}\label{lem:criticalmerit}
Let $\sigma_k=\sigma_k(\bar x)$. There is $\Acrit\ge0$ such that
\begin{equation}\label{eq:criticalmerit2}
 t_k^2(\sigma_k+q_k)\longrightarrow\Acrit,
 \qquad
 0\le U_k-\sigma_k=o(k^{-2}),
 \qquad
 t_k^2(U_k+q_k)\longrightarrow\Acrit.
\end{equation}
\end{lemma}
\begin{proof}
For fixed $z\in\mathcal C$, \eqref{eq:projection}, $x^k\to\bar x$, and $v^k=O(k^{-1})$ give
\begin{equation}\label{eq:energyasympt}
 \EE_k(z)=t_k^2(\sigma_k(z)+q_k)+\frac{\ell}{2}\norm{\bar x-z}^2+o(1).
\end{equation}
Indeed, $w_k-\bar x=(x^k-\bar x)+(t_k-1)v^k$, and its squared norm equals $t_k^2D_k+o(1)$; the cross term with $\bar x-z\in\mathcal V$ tends to zero by \eqref{eq:projection}. Thus
\[
 A(z):=\lim_k t_k^2(\sigma_k(z)+q_k)
\]
exists and is nonnegative. Write $\Acrit=A(\bar x)$. By \eqref{eq:activeC}, $A(z)$ is a common value $A_{\rm ri}$ for all $z\in\ri\mathcal C$, and $A(z)\le A_{\rm ri}$ for every $z\in\mathcal C$.

We prove $A_{\rm ri}=\Acrit$. Choose $z_0\in\ri\mathcal C$ and $z_\tau=(1-\tau)\bar x+\tau z_0$ for $0<\tau\le1$. These points belong to $\ri\mathcal C$. From the corrected monotonicity, for every fixed $K$,
\[
 A_{\rm ri}+\frac{\ell}{2}\norm{\bar x-z_\tau}^2
 =\lim_k\EE_k(z_\tau)\le\EE_K(z_\tau)+\epsilon_K.
\]
For each $i$, convexity along the segment and lower semicontinuity give $F_i(z_\tau)\to F_i(\bar x)$ as $\tau\downarrow0$, hence $\EE_K(z_\tau)\to\EE_K(\bar x)$. Letting $\tau\downarrow0$ and then $K\to\infty$ yields $A_{\rm ri}\le\Acrit$. The reverse inequality was already established. Therefore
\begin{equation}\label{eq:Azbound}
 A(z)\le\Acrit,
 \qquad z\in\mathcal C.
\end{equation}

It remains to control the moving comparators in the merit supremum. Use the bounded set
\[
 Z=\{z\in X:F(z)\in F(\Xstar\cap\mathcal L_0),\ \norm{z-x^0}^2\le R+1\}.
\]
For each $N$, choose $z_N\in Z$ with
\begin{equation}\label{eq:nearmax}
 0\le\sigma_N(z_N),
 \qquad
 U_N-N^{-3}\le\sigma_N(z_N)\le U_N.
\end{equation}
If $U_N=0$, a representative dominating $x^N$ gives such a point. If $U_N>0$, positive comparisons and their weak Pareto representatives give the required approximation. Thus $F(z_N)\le F(x^N)$.

Along any subsequence with $z_N\to z_*$, lower semicontinuity gives $z_*\in\mathcal C$. For $K<N$, the backward inequality and \eqref{eq:nearmax} imply $\sigma_j(z_N)\ge-q_j$ for $K\le j<N$. Summing \eqref{eq:energy-recursion} over this interval gives
\begin{equation}\label{eq:restartbound}
 \EE_N(z_N)\le\EE_K(z_N)+\epsilon_K.
\end{equation}
For fixed $K$, the map $z\mapsto\sigma_K(z)$ is upper semicontinuous, so $\limsup_N\EE_K(z_N)\le\EE_K(z_*)$. Boundedness of $w_N$ and \eqref{eq:projection} also give
\[
 \EE_N(z_N)-\EE_N(\bar x)
 =t_N^2(\sigma_N(z_N)-\sigma_N)
 +\frac{\ell}{2}\norm{\bar x-z_*}^2+o(1).
\]
Combining with \eqref{eq:restartbound} and sending $K\to\infty$ yields
\[
 \limsup_N t_N^2(\sigma_N(z_N)-\sigma_N)
 \le A(z_*)-\Acrit\le0.
\]
The scaled differences are bounded because $U_N=O(N^{-2})$, $\sigma_N\ge-q_N$, and $q_N=O(N^{-2})$. A subsequence realizing their limsup has a further convergent comparator subsequence because $Z$ is bounded in finite dimension, so the same conclusion holds for the full limsup. Since $U_N\ge\sigma_N$ and $t_N^2N^{-3}\to0$, \eqref{eq:nearmax} gives $t_N^2(U_N-\sigma_N)\to0$. The last assertion follows from the first.
\end{proof}

Next we introduce the defects $\hh_k$ and $\bb_k$ that measure the failure of the energy recursion to be exact. The defect $\omega_k$ measures the gap between the model and its linearization at $x^k$. The following proposition gives the budgets for these defects and their asymptotic behavior.

\begin{proposition}\label{prop:criticaldefects}
Under the critical rule,
\begin{equation}\label{eq:budgets}
 \hh_k\ge0,
 \quad
 \bb_k\ge0,
 \quad
 0\le\omega_k\le\hh_k,
 \qquad
 \sum_{k\ge1}k^2\hh_k<\infty,
 \qquad
 \sum_{k\ge1}k\bb_k<\infty,
\end{equation}
where, with $p_k=x^{k+1}$, $r_k=p_k-y^k$, $\xi_k=-\ell r_k$, $\psi_k(z)=\psi(z;x^k,y^k)$,
\begin{align}
 \hh_k&=\beta_k^2q_k-q_{k+1}-(\sigma_{k+1}-\sigma_k),\label{eq:defecth}\\
 \bb_k&=\frac{\ell}{2}\left(\norm{y^k-\bar x}^2-\norm{p_k-\bar x}^2\right)-\sigma_{k+1},\label{eq:defectb}\\
 \omega_k&=\ip{\xi_k}{p_k-x^k}-\psi_k(p_k).\label{eq:defectomega}
\end{align}
The symbols $\bb_k$ and $\hh_k$ denote defects, not momentum parameters. For every fixed $\rho>1$,
\begin{equation}\label{eq:windowlimits}
 \sum_{k=K}^{\lfloor\rho K\rfloor}kq_k\longrightarrow\frac{4\Acrit}{3}\log\rho,
 \qquad
 \sum_{k=K}^{\lfloor\rho K\rfloor}kU_k\longrightarrow\frac{8\Acrit}{3}\log\rho.
\end{equation}
In particular, $U_k=o(k^{-2})$ if and only if $\Acrit=0$; in that case also $\norm{v^k}=o(k^{-1})$.
\end{proposition}
\begin{proof}
The nonnegativity of $\hh_k$ and $\bb_k$ follows from \eqref{eq:gap-monotone} and \eqref{eq:distance}. By \eqref{eq:aggregate-opt}, $\xi_k\in\partial\psi_k(p_k)$. Since $\psi_k(x^k)\le0$, the subgradient inequality gives $\omega_k\ge-\psi_k(x^k)\ge0$. The identity
\[
 \beta_k^2q_k-q_{k+1}
 =\ip{\xi_k}{p_k-x^k}+\frac{\ell}{2}\norm{r_k}^2
\]
implies $\hh_k-\omega_k=Q_k(p_k)-(\sigma_{k+1}-\sigma_k)\ge0$.

Keeping the defects in the squared-norm calculation that leads to \eqref{eq:energy-recursion} gives the exact identity
\begin{equation}\label{eq:exactdefectenergy}
 \EE_k(\bar x)-\EE_{k+1}(\bar x)
 =\frac14\sigma_k+t_{k+1}(t_{k+1}-1)\hh_k+t_{k+1}\bb_k.
\end{equation}
The energies converge and $\sum_k|\sigma_k|<\infty$ by \cref{lem:criticalpoint}. Summing \eqref{eq:exactdefectenergy} proves the two budgets.

Let $V_k=\sigma_k+q_k\ge0$. Then
\begin{equation}\label{eq:Vrecurrence}
 V_k-V_{k+1}=\eta_kq_k+\hh_k,
 \qquad
 k^2V_k\longrightarrow4\Acrit.
\end{equation}
Put $M=\lfloor\rho K\rfloor$. Summation by parts gives
\begin{align*}
 \sum_{k=K}^M k^2(V_k-V_{k+1})
 &=K^2V_K-M^2V_{M+1}+\sum_{k=K+1}^M(2k-1)V_k\\
 &\longrightarrow8\Acrit\log\rho.
\end{align*}
The two boundary limits cancel. Moreover, $k^2\eta_k=6k+O(1)$, $\sum_kq_k<\infty$, and $\sum_kk^2\hh_k<\infty$. Substitution of \eqref{eq:Vrecurrence} gives the first limit in \eqref{eq:windowlimits}. Since $k^2(U_k+q_k)\to4\Acrit$ by \cref{lem:criticalmerit}, subtracting the first limit from the corresponding window sum for $U_k+q_k$ gives the second.

If $\Acrit=0$, nonnegativity and \eqref{eq:criticalmerit2} give both little-o rates. Conversely, $U_k=o(k^{-2})$ makes its fixed-ratio window sums tend to zero; the second limit in \eqref{eq:windowlimits} then forces $\Acrit=0$. The displacement conclusion follows from $0\le q_k\le V_k$.
\end{proof}

This proves \cref{thm:criticalenergy}. The critical little-o question is now exactly the question whether $\Acrit$ must vanish for general convex composite data; the equivalence separates this from the finite $O(N^{-2})$ baseline.

\section{Convergence rate for \texorpdfstring{$1<\alpha<3$}{1<alpha<3}}\label{sec:alpha-between}

For $1<\alpha<3$ the rational clock is positive and the signed energy gives a polynomial upper bound. This is the subcritical range: the root parameter $a=(\alpha-3)/(\alpha-1)$ in \eqref{eq:rational-map} is negative, so \cref{thm:root} does not apply, but the energy identity remains valid because it only uses the algebraic recurrence.

\begin{theorem}\label{thm:subcritical}
Under \cref{ass:geometry}, suppose the rational rule \eqref{eq:rational} is used with $1<\alpha<3$. Then
\begin{equation}\label{eq:subcritical}
 U_N=O\bigl(N^{-(\alpha-1)}\bigr).
\end{equation}
The estimate is an upper bound obtained from the signed energy; it is not claimed to be optimal.
\end{theorem}
\begin{proof}[Proof sketch]
Set $a=\alpha-1\in(0,2)$ and use the clock $t_k=(k+\alpha-2)/(\alpha-1)$. Define the signed energy with $w_k=t_kx^k-(t_k-1)x^{k-1}$ and the signed dissipation $\delta_k=t_k^2-t_{k+1}^2+t_{k+1}$. The basic estimates give a recursion of the form
\[
 \HH_{k+1}(z)\le\HH_k(z)-c_k\sigma_k(z)+\text{nonnegative displacement terms},
\]
where the coefficient $c_k$ is eventually nonnegative. Multiplying the scalarized version by a positive multiplier and summing produces a Gamma ratio: for some $K$ and all $N\ge K$,
\begin{equation}\label{eq:gamma}
 U_N\le M_K\frac{\Gamma(N)\Gamma(K+a)}{\Gamma(K)\Gamma(N+a)}
 =O(N^{-a}),
 \qquad a=\alpha-1.
\end{equation}
The constant $M_K$ depends on the fixed prefix and on the problem data. The argument gives only the stated upper bound; it neither proves optimality nor excludes a faster rate. At the endpoint $\alpha=3$, the finite critical bound $2\ell R/[N(N+3)]$ of \cref{thm:critical} is sharper than \eqref{eq:gamma}.
\end{proof}

In terms of Table~\ref{tab:rates}, the subcritical range has a rate that degrades continuously as $\alpha$ decreases to $1$, and no displacement rate is asserted. This is the multiobjective counterpart of the scalar subcritical behaviour, where the exponent also depends linearly on the damping parameter.

\section{Convergence for \texorpdfstring{$\alpha>0$}{alpha>0}}\label{sec:alpha-positive}

The following qualitative conclusion holds for the whole positive range. It does not assert a power rate or point convergence.

\begin{theorem}[Qualitative convergence]\label{thm:qualitative}
Under \cref{ass:geometry}, for every $\alpha>0$ the rational rule \eqref{eq:rational} satisfies
\begin{equation}\label{eq:qualitative}
 U_k\longrightarrow0,
 \qquad
 \norm{v^k}\longrightarrow0,
 \qquad
 F_i(x^k)\longrightarrow c_i\in\R
 \quad(i=1,\ldots,m),
\end{equation}
and
\begin{equation}\label{eq:sumDk}
 \sum_{k\ge1}\frac{D_k}{k}<\infty.
\end{equation}
\end{theorem}
\begin{proof}[Proof sketch]
Let $e_{i,k}=F_i(x^k)+q_k$. The componentwise inequality \eqref{eq:component} and $\beta_k^2\le1$ show that $e_{i,k}$ is nonincreasing. The trajectory remains in the initial lower level set and is bounded below on that set by properness and lower semicontinuity, so each $e_{i,k}$ has a finite limit; since $q_k\to0$, also $F_i(x^k)\to c_i$.

For the rational rule,
\[
 1-\beta_k^2=1-\left(\frac{k-1}{k+\alpha-1}\right)^2
 =\frac{2\alpha k+\alpha^2-2\alpha}{(k+\alpha-1)^2}
 \asymp\frac{2\alpha}{k}.
\]
Multiplying the signed componentwise inequality by a scaling factor and summing, with these coefficients, gives \eqref{eq:sumDk}. The same scaling makes $F_i(x^k)+q_k$ monotone and bounded below, so it has a finite limit. A tail argument that keeps the signed comparison terms then forces $U_k\to0$ and $D_k\to0$. The argument does not require boundedness of the initial lower level set.
\end{proof}

The critical point convergence in \cref{thm:criticalenergy} is a separate result. It is not a consequence of \cref{thm:qualitative}: the qualitative theorem gives $U_k\to0$ and $D_k\to0$ but not convergence of the full sequence or a critical energy representation.

\section{Numerical illustration of the little-\texorpdfstring{$o$}{o} rate}\label{sec:numerical}

We illustrate the distinction between the $O(k^{-2})$ baseline and the $o(k^{-2})$ rate of \cref{thm:root} on the composite biobjective JOS1 instance used in the benchmark of \cite{tanabe2023accelerated}:
\begin{equation}\label{eq:jos1}
 F_1(x)=f_1(x)+g_1(x),
 \qquad
 F_2(x)=f_2(x)+g_2(x),
\end{equation}
with
\begin{equation}\label{eq:jos1-parts}
 f_1(x)=\frac{\norm{x}^2}{n},
 \quad
 f_2(x)=\frac{\norm{x-2e}^2}{n},
 \quad
 g_1(x)=\frac{\norm{x}_1}{n},
 \quad
 g_2(x)=\frac{2\norm{x-e}_1}{n},
\end{equation}
where $e=(1,\ldots,1)^\top$ and $n=1000$. The smooth parts have Lipschitz gradients with $L=2/n=2\cdot10^{-3}$; the nonsmooth parts are the shifted $\ell_1$ terms of the JOS1 benchmark.

\begin{figure}[!b]
\centering
\includegraphics[width=\textwidth]{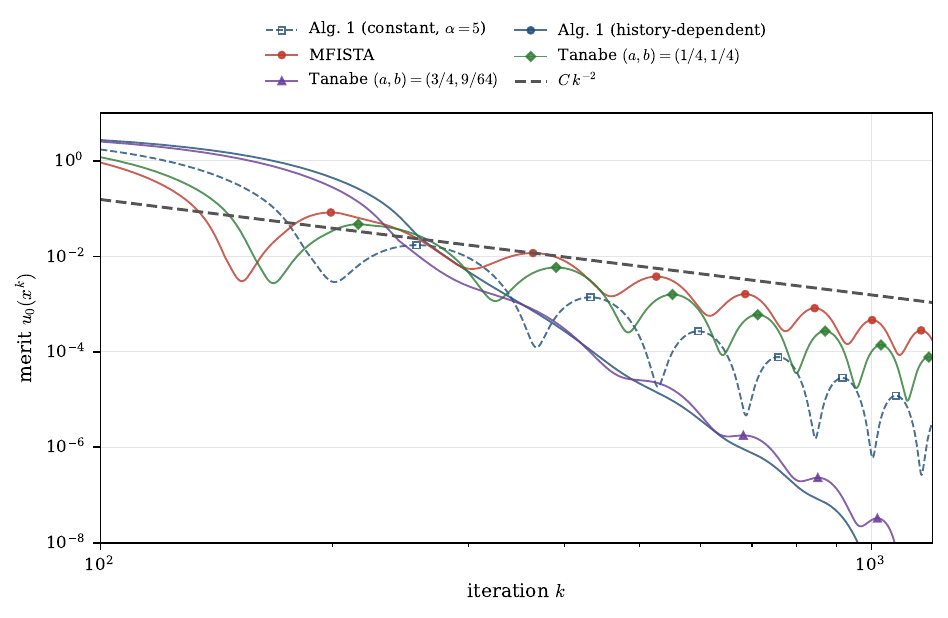}
\caption{Composite JOS1 instance \eqref{eq:jos1}--\eqref{eq:jos1-parts} with $n=1000$ and $\ell=5$. The full trajectories $k\mapsto u_0(x^k)$ are drawn as solid and dashed curves; the markers indicate their local maxima. The dashed line is $Ck^{-2}$, where $C=1.56\times10^{3}$ is the smallest constant for which the line is an upper envelope of the MFISTA maxima for $k\ge200$. The constant $\alpha=5$ instance, the history-dependent rule \eqref{eq:history-rule}, and the two Tanabe instances all move below the reference, while the scaled merit of MFISTA remains at the level of the reference. The history-dependent rule tracks the tuned $(3/4,9/64)$ pair and is competitive with it, while the $(1/4,1/4)$ pair lies between the Algorithm~1 curves and MFISTA; all of them fall below the reference, consistently with the $o(k^{-2})$ bound of \cref{thm:root}.}
\label{fig:jos1-rate}
\end{figure}

The merit function is evaluated exactly. The data are invariant under permutations of the coordinates, and every weighted-sum minimizer is of the form $ce$. Hence every weak Pareto point is of this form, and along $ce$ the two objectives are
\begin{equation}\label{eq:jos1-front}
 F_1(ce)=c^2+|c|,
 \qquad
 F_2(ce)=(c-2)^2+2|c-1|.
\end{equation}
Consequently
\begin{equation}\label{eq:jos1-merit}
 u_0(x)=\max_{c\in\R}\min_{i=1,2}\bigl(F_i(x)-F_i(ce)\bigr).
\end{equation}
The function $c\mapsto\min_i(F_i(x)-F_i(ce))$ is concave, so \eqref{eq:jos1-merit} is a one-dimensional concave maximization, which we evaluate exactly by a golden-section search.

For \cref{alg:root} we use the history-dependent rule
\begin{equation}\label{eq:history-rule}
\begin{aligned}
 \phi_k&=\frac{D_k}{D_k+D_{k-1}},
 &\qquad
 \theta_k&=\frac{D_k}{D_k+S_k},\\
 a_k&=\frac34+\frac15\,\phi_k\,(1+k)^{-1/4},
 &\qquad
 b_k&=\frac{a_k^2}{4}+\left(\frac14-\frac{a_k^2}{4}\right)\theta_k,
\end{aligned}
\end{equation}
with $D_k=\norm{x^k-x^{k-1}}^2$, $S_k=\sum_{j=1}^kD_j$, and $D_0=D_1$. Since $0\le\phi_k,\theta_k\le1$, we have $b_k\in[a_k^2/4,1/4]$ for every $k$; moreover $\theta_k\to0$, hence $a_k\to3/4$ and $b_k\to a_k^2/4=9/64$. Thus the history-dependent rule adapts towards the tuned Tanabe pair $(3/4,9/64)$ used below. It uses only the displacements already computed at step $k$, so it is a genuine history-dependent choice in the sense of \cref{alg:root}.

We compare five momentum rules with the same step parameter $\ell=5$ and the same initial point, drawn uniformly from $[-2,4]^n$ with seed $2026$:
\begin{itemize}[leftmargin=2.2em]
\item the constant rational instance $\alpha=5$ of \eqref{eq:rational-map}, that is, $a_k\equiv1/2$ and $b_k\equiv1/16$;
\item the history-dependent rule \eqref{eq:history-rule};
\item MFISTA of \cite{tanabe2023accelerated}, the classical accelerated rule $t_{k+1}=\sqrt{t_k^2+1/4}+1/2$, which is the $O(k^{-2})$ baseline;
\item the constant root instance $(a_k,b_k)\equiv(1/4,1/4)$ from Table~1 of \cite{tanabe2022root};
\item the constant root instance $(a_k,b_k)\equiv(3/4,9/64)$, the best JOS1-L1 row of Table~1 in \cite{tanabe2022root}.
\end{itemize}
The first, fourth, and fifth rules are constant-parameter instances of \cite{tanabe2022root}; the second is the history-dependent rule analyzed in this paper.
The choice $\ell=5$ is deliberately conservative. It satisfies $\ell\ge L$ and keeps the asymptotic regime visible over the displayed interval; the exponent in the rate statements does not depend on the size of $\ell$. All subproblems are solved through their one-dimensional dual, as in Section~6 of \cite{tanabe2023accelerated}.

Accelerated iterates are non-monotone, so the full trajectories oscillate between values close to the Pareto front and overshoots away from it. The local maxima are the quantities controlled by a worst-case rate, and they are the markers in \cref{fig:jos1-rate}. Over the displayed maxima the scaled merit of the constant $\alpha=5$ instance decreases from $k^2u_0(x^k)\approx1.1\times10^{3}$ at $k=257$ to $\approx1.4\times10^{1}$ at $k=1076$, a factor of about $80$. For MFISTA they are $1.6\times10^{3}$ at $k=364$ and $3.8\times10^{2}$ at $k=1160$, a factor of about $4$; for the Tanabe $(1/4,1/4)$ instance they are $2.2\times10^{3}$ at $k=216$ and $1.1\times10^{2}$ at $k=1187$, a factor of about $20$. The history-dependent rule is slower in the early transient ($u_0\approx4.4\times10^{-1}$ at $k=200$) but then decreases steeply to $2.7\times10^{-4}$ at $k=400$, $3.7\times10^{-6}$ at $k=600$, and $1.5\times10^{-7}$ at $k=800$, tracking the tuned Tanabe pair. The tuned Tanabe $(3/4,9/64)$ instance is slower during the initial transient ($u_0\approx2.8\times10^{-1}$ at $k=200$) but then decreases steeply to $3.0\times10^{-4}$ at $k=400$ and $2.0\times10^{-7}$ at $k=800$, and, together with the history-dependent rule, it is the fastest at the end of the range; this is consistent with it being the best JOS1-L1 row of Table~1 in \cite{tanabe2022root}. Thus all five rules move below the $k^{-2}$ reference at the end of the range, while the scaled merit of MFISTA remains at the level of the reference. The history-dependent rule is tuned to the same pair as the tuned Tanabe instance; its trajectory confirms that the little-$o$ guarantee is compatible with a competitive history-dependent choice. The generated parameters satisfy $a_k\in[0.767,0.834]$ and $b_k\in[0.147,0.212]$, with $a_k$ decreasing toward $3/4$ and $b_k$ decreasing toward $a_k^2/4=9/64$. The experiment illustrates the asymptotic result rather than replacing its proof in \cref{thm:root}.

\section{Conclusion}\label{sec:conclusion}

We have developed a unified analysis of accelerated proximal gradient methods with general inertial term for convex composite multiobjective optimization. The analysis keeps the original offset subproblem, the original composite splitting, and the step size $1/\ell$ with $\ell\ge L$. The rate picture is organized by the momentum coefficient, in the same way as the scalar picture organized by the damping parameter, and it is summarized in Table~\ref{tab:rates}.

For history-dependent root momentum, the chain of basic estimates, energy identity, coefficient asymptotics, and auxiliary lemmas yields boundedness, weighted square summability of the displacements, $U_k=o(k^{-2})$, and $\norm{v^k}=o(k^{-1})$ whenever $0\le a_k<1$, $a_k^2/4\le b_k\le1/4$, and $a_k\to\bar a\in(0,1)$. The second parameter need not converge or have bounded total variation. A finite negative prefix of the dissipation coefficient is retained as an explicit correction and does not restart the method.

At the critical parameter $\alpha=3$, the same machinery gives the finite merit bound $2\ell R/[N(N+3)]$, square summability of the displacements, displacement decay $O(1/N)$, convergence of the full sequence to a weak Pareto point, and a critical energy $\Acrit$ such that $U_k=o(k^{-2})$ exactly when $\Acrit=0$. The remaining critical little-o question is therefore a question about the vanishing of a nonnegative limit, not a missing large-$O$ estimate.

For the subcritical range, $1<\alpha<3$ has the rate $O(N^{-(\alpha-1)})$, and every $\alpha>0$ has qualitative convergence with $\sum_kD_k/k<\infty$. The resulting exponents decrease continuously with $\alpha$, which is the vector analogue of the subcritical behaviour of the accelerated gradient method.

For the future, we will try to find how to extend the results to the nonconvex case, and to study the convergence of the whole sequence for $\alpha>3$. And we want to design a unified framework to construct the auxiliary functions for the accelerated multiobjective methods and even other multiobjective methods.







\begingroup
\renewcommand{\bibfont}{\small}

\endgroup

\end{document}